\documentclass{elsart}
\usepackage{amssymb}
\newtheorem{lemma}{Lemma}[section]
\newtheorem{theorem}{Theorem}[section]

\begin{document}
\begin{frontmatter}
\title{Lyapunov-Sylvester operators for Generalized Nonlinear Euler-Poisson-Darboux System}
\author{Anouar Ben Mabrouk}
 \ead{anouar.benmabrouk@fsm.rnu.tn}
% \corauth[cor1]{Corresponding author}
 \address{UR Algebra, NUmber Theory and Nonlinear Analysis, Department of Mathematics, Faculty of Sciences, 5019 Monastir.Tunisia}
%%%%%%%%%%%%%%%%%%%%%%%%%%%%%%%%%%%%%%%%%%%%%%%%%%%%%%%%%%%%%%%%%%%%%%%%%%%%%%%%%%%%%%%%%%%%%%%%%%%%%%%%%%%%%%%%%%%%%%%%%
%%%%%%%%%%%%%%%%%%%%%%%%%%%%%%%%%%%%%%%%%%%%%%%%%%%%%%%%%%%%%%%%%%%%%%%%%%%%%%%%%%%%%%%%%%%%%%%%%%%%%%%%%%%%%%%%%%%%%%%%%
%%%%%%%%%%%%%%%%%%%%%%%%%%%%%%%%%%%%%%%%%%%%%%%%%%%%%%%%%%%%%%%%%%%%%%%%%%%%%%%%%%%%%%%%%%%%%%%%%%%%%%%%%%%%%%%%%%%%%%%%%
\begin{abstract}
In this paper a nonlinear Euler-Poisson-Darboux system is considered. In a first part, we proved the genericity of the hypergeometric functions in the development of exact solutions for such a system in some special cases leading to Bessel type differential equations. Next, a finite difference schemle in two-dimensional case has been developed. The continuous system is transformed into an algebraic quasi linear discrete one leading to generalized  Lyapunov-Sylvester operators. The discrete algebraic system is proved to be uniquely solvable, stable and convergent based on Lyapunov criterion of stability and Lax-Richtmyer equivalence theorem for the convergence. A numerical example has been provided at the end to illustrate the efficiency of the theoretical results. The present method is thus proved to be more accurate than existing ones and lead to faster algorithms.
\end{abstract}
\end{frontmatter}
\section{Introduction}
In this work we use Lyapunov-Sylvester algebraic operators to approximate the solutions of a generalized Euler-Poisson-Darboux (EPD) system in two-dimensional case. The present article is precisely devoted to the development of a numerical method based on two-dimensional finite difference scheme to approximate the solution of the generalized EPD system in $\mathbb{R}^2$ in the presence of mixed power laws nonlinearities. Denote for $a\in\mathbb{R}$, $\Gamma_a(x)=\displaystyle\frac{2a}{x}$ and for $\lambda,\gamma$ in $\mathbb{R}$, $F_{\lambda,\gamma}(x)=(\Gamma_\lambda(x),\Gamma_\gamma(x))$. We consider the following evolutive system.
\begin{equation}\label{ContinuousProblem1}
\left\{\begin{array}{lll}
u_{tt}+\Gamma_a(t)v_t=\Delta u+<F_{\lambda,\gamma},\nabla v>+|u|^{p-1}v,\\
v_{tt}+\Gamma_a(t)u_t=\Delta v+<F_{\lambda,\gamma},\nabla u>+|v|^{q-1}u
\end{array}\right.
\end{equation}
with initial conditions
\begin{equation}\label{eqn1-3}
u(x,y,t_{0})=u_{0}(x,y)\quad \hbox{and}\quad\frac{\partial u}{\partial t}(x,y,t_{0})=u_1(x,y),\quad(x,y)\in\Omega
\end{equation}
and boundary conditions
\begin{equation}\label{eqn1-4}
\frac{\partial u}{\partial\eta}(x,y,t)=0,\quad((x,y),t)\in\partial\Omega\times(t_{0},+\infty).
\end{equation}
on a rectangular domain $\Omega =[L_{0},L_1]\times[L_{0},L_1]$ in $\mathbb{R}^2$. $t_{0}\geq 0$ is a real parameter fixed as the initial time, $u_{t}$ is the first order partial derivative in time, $u_{tt}$ is the second order partial derivative in time, $\Delta =\frac{\partial ^2}{\partial x^2}+\frac{\partial ^2}{\partial y^2}$ is the Laplace operator on $\mathbb{R}^2$. $\frac{\partial}{\partial\eta}$ is the outward normal derivative operator along the boundary $\partial \Omega $. Finally, $u$, $u_{0}$ and $u_1$ are real valued functions with $u_{0}$ and $u_1$ are $\mathcal{C}^2$ on $\overline{\Omega}$. $u$ and $v$ are the unknown candidates supposed to be $\mathcal{C}^{4}$ on $\overline{\Omega}$. $p$ and $q$ are reel parameters such that $p,q>1$.

Several papers have been devoted to the study of existence and uniqueness of solutions of the linear problem \ref{ContinuousProblem1} (without the nonlinear parts) for the case where $u\equiv v$ and where the system is in fact the well known EPD equation. In some of these studies, exact solutions are developed such as solitary, stationary, time-independent, one-dimensional ones. For example, in the case of a one-direction viscous fluid we may seek solutions of the form $u(x,y,t)=\alpha\psi(x)$. In this case, the problem is transformed into a one variable ordinary differential equation
$$
q\psi''(x)+\psi(x)+\alpha\psi^2(x)=ax+b,
$$
for some constants $a$ and $b$ depending on the initial-boundary conditions. The existence and uniqueness problems are overcame using ODEs. 

In \cite{Seilkhanova1}, a generalized form of EPD equation is considered and particular solutions are constructed in an explicit form expressed by the Lauricella hypergeometric function of three variables. Properties of each constructed solutions have been investigated in sections of surfaces of the characteristic cone. The authors proved that found solutions have singularity $\frac{1}{r}$ as $r\rightarrow0$, where $r^2=(x-x_0)^2+(y-y_0)^2-(t-t_0)^2$.

In \cite{Urinov}, a singular Cauchy problem for the multi-dimensional EPD equation with spectral parameter has been investigated with the help of the generalized Erdelyi-Kober fractional operator. Solution of the considered problem is found in explicit form for various values of the parameter $p$ of the equation.

In the present work, we intend to apply some algebraic operators to develop numerical solutions for EPD system. It consists of the well known Lyapunov-Sylvester operators. For given matrices $A\in \mathbb{R}^{m\times m},B\in \mathbb{R}^{n\times n}$, and $C\in \mathbb{R}^{m\times n}$, the Sylvester equation is given by the form $AX+XB=C$. A classical idea to obtain the solution $X$ is to rewrite the Sylvester equation in standard\ $mn\times mn$ linear system $G\widetilde{x}=\widetilde{c}$ using the Kronecker Product,\cite{Jameson}. The Sylvester equation can be solved by Gaussian elimination with $O(m^{3}n^{3})$ flops. This approach dramatically increases the complexity of the computation, and also cannot preserve the intrinsic properties of the problem in practice \cite{Simoncini}. 

In numerical analysis, for solving the Sylvester equation one using the Bratels-Stewart and the Golub-Nash-Van Loan algorithm use $O(m^{3}+n^{3})$ floating point operations, if one assume that an $m\times m$ matrix can be reduced to Schur form with $O(m^{3})$ operations (See \cite{Bartels} and \cite{Goloub}).

In \cite{Kirrini} the author described an algorithm that computes the solution $X$ over an arbitrary field $\mathbb{F}$. The complexity of the algorithm for $A\in \mathbb{F}^{m\times m},\ \ B\in\mathbb{F}^{n\times n}$ and $m,n\leq N$ is $O(N^{\beta }\log N)$ arithmetic operations in $\mathbb{F}$, where $\beta>2$ is such that $M\times M$ matrices can be multiplied with $O(M^{\beta}) $ arithmetic operations. This algorithm is competitive in terms of arithmetic operation with and even faster than the classical algorithms.

The method developed in this paper consists in replacing time and space partial derivatives by finite-difference approximations in order to transform the continuous problem into quasi linear Lyapunov-Sylvester system. The motivation behind the application of Lyapunov-Syslvester operators was already evoked in \cite{Benmabrouk1}. We recall in brief that such a method leads to fast convergent and more accurate discrete algebraic systems without going back to the use of tri-diagonal and/or fringe-tridiagonal matrices already used when dealing with multidimensional problems especially in discrete PDEs.

To recapitulate, the method developed here is favorable for many reasons.
\begin{itemize}
	\item The first motivation is the fact that it somehow does not change the geometric presentation of the problem as we propose to solve in the same two-dimensional space. We did not project the problem on tri-diagonal representations using the Kronecker product. Relatively to computer architecture, the process of projecting on different spaces and next lifting to the original one may induce degradation of error estimates and slow algorithms.
	\item The method developed is not just a resolution of a PDE. We recall already that the resolution itself is not a negligible aim. Further, it proves the efficiency of algebraic operators other than classical tri-diagonal ones.
	\item We proved here that even when the two systems are equivalent in the sense that they present the same PDE, but with different forms and dimensions, such forms play a major role in the resolution.
	\item The fact of obtaining fast algorithms is very important in computer sciences and makes itself a major aim in computer studies. Recall that the famous method known in mathematical studies of accelerating algorithms in the EM one (expectation-maximisation) which is based on more complicated theories. Here, we proved that we may obtain more rapid algorithms by using just a suitable representation and suitable discerete transformation of the PDE. We got faster algorithms without adding more parameters.
\end{itemize}
In the present work, existence and multiplicity of the solutions of problem (\ref{ContinuousProblem1})-(\ref{eqn1-4}) are developed in some special cases. We showed that special functions such as hypergeometric series and Bessel function are generic for developing such solutions. Section 3 is devoted to the development of a 2-dimensional discrete scheme to transform the continuous problem (\ref{ContinuousProblem1})-(\ref{eqn1-4}) to a discrete one. A system of generalized Lyapunov-Sylvester equations is obtained. The solvability of such a discrete system is proved next in section 4. Section 5 is concerned with the consistency, stability and the convergence of the discrete Lyapunov-Sylvester problem obtained in section 3. The crucial idea is the application of the truncation error for consistency, Lyapunov cretirion for stability and the Lax equivalence theorem for the convergence. Section 6 is devoted to the development of a numerical example. The performance of the discrete scheme is proved by means of error estimates as well as fast algorithms. The conclusion is finally subject of section 7. 

\section{On the existence and multiplicity of the solutions of the continuous problem}
In this part, we review the generalized linear EPD system. We will show the genericity of hypergeometric functions in the development of solutions. The generalized linear EPD system is
\begin{equation}\label{LinearContinuousProblem}
\left\{\begin{array}{lll}
u_{tt}+\Gamma_a(t)v_t=\Delta u+<F_{\lambda,\gamma},\nabla v>,\\
v_{tt}+\Gamma_a(t)u_t=\Delta v+<F_{\lambda,\gamma},\nabla u>.
\end{array}\right.
\end{equation}
\subsection{A first class of time-independent solutions}
In this section we propose to develop a first class of solutions of problem (\ref{LinearContinuousProblem}) on the time-independent additive form
$$
u(x,y,t)=v(x,y,t)=\varphi(x,y)=f(x)+g(y).
$$
The stationary problem associated to the system (\ref{LinearContinuousProblem}) is
\begin{equation}\label{ContinuousProblem2}
\Delta u+<F_{\lambda,\gamma},\nabla v>=\Delta v+<F_{\lambda,\gamma},\nabla u>=0.
\end{equation}
Substituting $f$ and $g$ in the last equations yield that
$$
f''(x)+\frac{2\lambda}{x}f'(x)=g''(y)+\frac{2\gamma}{y}g'(y).
$$
Hence, there exists a constant $K\in\mathbb{C}$ for which
\begin{equation}\label{ContinuousProblem3}
\left\{\begin{array}{lll}
f''(x)+\frac{2\lambda}{x}f'(x)=K\\
g''(y)+\frac{2\gamma}{y}g'(y)=-K.
\end{array}\right.
\end{equation}
Therefore, whenever $\lambda,\gamma\notin\{\pm\frac{1}{2}\}$, by applying classical resolution of ODEs we get the general solutions
$$
f(x)=\frac{K_1}{1-2\lambda}|x|^{1-2\lambda}sign(x)+\frac{K}{1+2\lambda}\frac{x^2}{2},\;x\not=0
$$
and
$$
g(y)=\frac{K_2}{1-2\gamma}|y|^{1-2\gamma}sign(y)-\frac{K}{1+2\gamma}\frac{y^2}{2},\;y\not=0.
$$
For the case when $\lambda=\gamma=\displaystyle\frac{1}{2}$ we get
$$
f(x)=K_1\log|x|+\frac{K}{6}x^3,\;x\not=0
$$
and
$$
g(y)=K_2\log|y|+\frac{K}{6}y^3,\;y\not=0.
$$
Similarly, for the case $\lambda=\gamma=-\displaystyle\frac{1}{2}$ we get
$$
f(x)=\frac{x^2}{2}\left(K\log|x|+K_1-\frac{1}{2}\right),\;x\not=0
$$
and
$$
g(y)=\frac{y^2}{2}\left(-K\log|y|+K_2-\frac{1}{2}\right),\;y\not=0.
$$
The remaining solutions (($\lambda\not=\pm\frac{1}{2}$ and $\gamma=\pm\frac{1}{2}$) and ($\lambda=\pm\frac{1}{2}$ and $\gamma\not=\pm\frac{1}{2}$)) are composed of linear combinaisions of the developed cases above. 
\subsection{Second class of time-independent solutions}
In this section, we continue to develop a second class of time-independent but multiplicative solutions. We consider solutions of the form
\begin{equation}\label{secondsolutionform}
u(x,y,t)=v(x,y,t)=\varphi(x,y)=f(x)g(y).
\end{equation}
Substituting as for the previous section $\varphi$ in the stationary system (\ref{LinearContinuousProblem}) we get
$$
g(y)\left(f''(x)+\frac{2\lambda}{x}f'(x)\right)=f(x)\left(g''(y)+\frac{2\gamma}{y}g'(y)\right).
$$
Hence, there exists a constant $K\in\mathbb{C}$ for which
\begin{equation}\label{ContinuousProblem4}
\left\{\begin{array}{lll}
f''(x)+\frac{2\lambda}{x}f'(x)=Kf(x)\\
g''(y)+\frac{2\gamma}{y}g'(y)=-Kg(y).
\end{array}\right.
\end{equation}
For example, when $\lambda=\gamma=\pm\displaystyle\frac{1}{\sqrt2}$, an explicit solution may be obtained by
\begin{equation}\label{solutionsinusoidale1}
f(x)=|x|^{-\lambda}\left(a_0\cos(\sqrt{K}x)+\frac{a_1}{\sqrt{K}}\sin(\sqrt{K}x)\right)
\end{equation}
and
\begin{equation}\label{solutionsinusoidale2}
g(y)=|y|^{-\gamma}\left(b_0\cos(\sqrt{K}y)+\frac{b_1}{\sqrt{K}}\sin(\sqrt{K}y)\right.
\end{equation}
As previously, the remaining solutions (($\lambda\not=\pm\frac{1}{\sqrt2}$ and $\gamma=\pm\frac{1}{\sqrt2}$) and ($\lambda=\pm\frac{1}{\sqrt2}$ and $\gamma\not=\pm\frac{1}{\sqrt2}$)) are composed of linear combinaisions of the solutions developed above. 
\subsection{A hypergeometric/Bessel type solution}
In this section we will prove that the EPD time-independent system has solutions that may be expressed by means of the famous special functions such as the hypergeometric one and its general variants especially Bessel function. For this we assume that the solutions $u$ and $v$ are on the form (\ref{secondsolutionform}) and that $f$ and $g$ are of the form
$$
y(x)=|x|^\nu\displaystyle\sum_{n=0}^{+\infty}a_nx^n.
$$
Substituting in the first equation of system (\ref{ContinuousProblem3}), we obtain the following recurrence system
\begin{equation}\label{recurence1}
\left\{\begin{array}{lll}
a_0(1+2\lambda\nu)=0,\\
a_1(1+2(\lambda+\nu)+2\lambda\nu)=0,\\
(1+2\lambda\nu+2(\lambda+\nu)n+n(n-1))a_n=Ka_{n-2},\;\forall\,n\geq2.
\end{array}\right.
\end{equation}
So, for $\lambda=-\nu=\pm\displaystyle\frac{1}{\sqrt2}$, we get
$$
a_n=\displaystyle\frac{K}{n(n-1)}a_{n-2},\;n\geq2
$$
which yields the same oscillating singular solutions provided in (\ref{solutionsinusoidale1}) and (\ref{solutionsinusoidale2}). \\
For $\lambda=-\nu$ but $1+2\lambda\nu\not=0$, we get firstly $a_0=a_1=0$. It remains in (\ref{recurence1}) just one equation to handle,
\begin{equation}\label{recurence2}
(1+2\lambda\nu+n(n-1))a_n=Ka_{n-2},\;\forall\,n\geq2.
\end{equation}
Assume now that there exists $p=2k\in\mathbb{N}$ for which $1+2\lambda\nu=-p(p-1)$. Equation (\ref{recurence2}) becomes
\begin{equation}\label{recurence3}
(n-p)(n+p-1)a_n=Ka_{n-2},\;\forall\,n\geq2
\end{equation}
which yields that $a_{2n+1}=0$, for all $n$, and that
$$
a_{2n}=\displaystyle\frac{1}{(p-1)!}\displaystyle\frac{(n+k)!K^{n-k}}{(n-k)!(2m+p)!}a_p
$$
and thus $f$ (and similarly $g$) may be expressed by means of a hypergeometric series
\begin{equation}\label{hypergeometricsolution1}
f(x)=\displaystyle\frac{a_p|x|^{-\lambda}}{(p-1)!}\displaystyle\sum_{n=k}^{+\infty}\displaystyle\frac{(n+k)!K^{n-k}}{(n-k)!(2m+p)!}x^{2n}.
\end{equation}
For $p=(2k+1)\in\mathbb{N}$ satisfying the same hypothesis $1+2\lambda\nu=-p(p-1)$, we get $a_{2n}=0$, for all $n$, and
$$
a_{2n+1}=(2p-1)!\displaystyle\frac{(n+k)!K^{n-k}}{(n-k)!(2m+p)!}a_p
$$
and thus $f$ (and similarly $g$) may be expressed by 
\begin{equation}\label{hypergeometricsolution2}
f(x)=a_p(2p-1)!|x|^{-\lambda}\displaystyle\sum_{n=k}^{+\infty}\displaystyle\frac{(n+k)!K^{n-k}}{(n-k)!(2m+p)!}x^{2n+1}.
\end{equation}
Assume now that $\lambda\not=\nu$ and that already $1+2\lambda\nu=0$, and denote $r=\lambda+\nu$ and $s=2r-1$. From equation (\ref{recurence1}) we obtain immediately $a_1=0$ and that
\begin{equation}\label{recurence4}
n(n+s)a_n=Ka_{n-2},\;\forall\,n\geq2.
\end{equation}
So, when $s\in-\mathbb{N}$, this yields that $a_{2n+}=0$, for all $n$ and that
$$
a_{2n}=\displaystyle\frac{\Gamma(r+1)}{\Gamma(s+1)}\displaystyle\frac{\Gamma(2n+s+2)K^{n}}{2^n\Gamma(n+r+1)}a_0.
$$
As a result,
\begin{equation}\label{hypergeometricsolution3}
f(x)=a_0\displaystyle\frac{\Gamma(r+1)}{\Gamma(s+1)}|x|^{-1/2\lambda}\displaystyle\sum_{n=k}^{+\infty}
\displaystyle\frac{\Gamma(2n+s+2)K^{n}}{2^n\Gamma(n+r+1)}x^{2n}.
\end{equation}
Now, when $\lambda\not=\nu$, $1+2\lambda\nu=0$ and $s=2(\lambda+\nu)-1=2r-1=-2k\in-\mathbb{N}$, the recurence relation (\ref{recurence1}) permits to obtain
$a_1=0$ and
\begin{equation}\label{recurence5}
n(n-p))a_n=Ka_{n-2},\;\forall\,n\geq2
\end{equation}
which in turn yields that $a_{2n+1}=0$, for all $n$ and that
$$
a_{2n}=a_{2k}k!\displaystyle\left(\displaystyle\frac{K}{4}\right)^{n-k}\displaystyle\frac{1}{(n-1)!(n-k)!}.
$$
As a result,
\begin{equation}\label{hypergeometricsolution4}
f(x)=a_{2k}k!\displaystyle\left(\displaystyle\frac{4}{K}\right)^{k}\displaystyle\left(\displaystyle\frac{\sqrt{K}|x|}{2}\right)^{-1/2\lambda}
\displaystyle\sum_{n=k}^{+\infty}\displaystyle\frac{1}{(n-1)!(n-k)!}x^{2n}.
\end{equation}
Next, for $\lambda\not=-\nu$, $1+2\lambda\nu\not=0$ but $1+2\lambda\nu+2(\lambda+\nu)=0$ and $\lambda+\nu=-k\in-\mathbb{N}$, the recurrence (\ref{recurence1}) becomes
\begin{equation}\label{recurence6}
(n-1)(n-2k)a_n=Ka_{n-2},\;\forall\,n\geq2.
\end{equation}
Therefore, $a_{2n}=0$ for $n\leq k-1$ and for $n\geq k$, we get
$$
a_{2n}=a_{2k}(2k)!\displaystyle\left(\displaystyle\frac{K}{2}\right)^{n-k}\displaystyle\frac{1}{(2n)!(n-k)!}.
$$
Similarly, for $n\geq k$ we get
$$
a_{2n+1}=a_{2k-1}(k-1)!\displaystyle\frac{K^{n-k}(n-k)!}{n!(2n-2k+1)!}.
$$
and for $n\leq k-1$, we have
$$
a_{2n+1}=\displaystyle\frac{a_{1}}{(2k-2)!}\displaystyle\left(\displaystyle\frac{-K}{2}\right)^{n}\displaystyle\frac{(2k-2n)!}{(k-n)!}.
$$
Thus we get for $\lambda\not=-1$,
\begin{equation}\label{hypergeometricsolution5}
\begin{array}{lll}
f(x)&=&\displaystyle\frac{a_{1}}{(2k-2)!}|x|^{-(1+2\lambda)/(2+2\lambda)}\displaystyle\sum_{n=0}^{k-1}\displaystyle\left(\displaystyle\frac{-K}{2}\right)^{n}
\displaystyle\frac{(2k-2n)!}{(k-n)!}x^{2n+1}\\
&&+a_{2k-1}(k-1)!|x|^{-(1+2\lambda)/(2+2\lambda)}\displaystyle\sum_{n=k}^{+\infty}\displaystyle\frac{K^{n-k}(n-k)!}{n!(2n-2k+1)!}x^{2n+1}\\
&&+a_{2k}(2k)!|x|^{-(1+2\lambda)/(2+2\lambda)}\displaystyle\sum_{n=k}^{+\infty}\displaystyle\left(\displaystyle\frac{K}{2}\right)^{n-k}
\displaystyle\frac{1}{(2n)!(n-k)!}x^{2n+1}.
\end{array}
\end{equation}
For $\lambda=-1$ and $\nu=\displaystyle\frac{1}{2}$, we get $a_0=0$ and
$$
a_{2n}=\displaystyle\left(\displaystyle\frac{K}{4}\right)^{n-1}\displaystyle\frac{1}{n!(n-1)!}a_2,\;n\geq1
$$
and
$$
a_{2n+1}=\displaystyle\frac{1}{K}\displaystyle\frac{(4K)^n(n!)^2}{(2n)!(2n+1)!}a_1,\;n\geq1.
$$
Hence,
\begin{equation}\label{hypergeometricsolution6}
\begin{array}{lll}
f(x)&=&a_2\sqrt{|x|}\displaystyle\sum_{n=1}^{+\infty}\displaystyle\left(\displaystyle\frac{K}{4}\right)^{n-1}\displaystyle\frac{1}{n!(n-1)!}x^{2n}\\
&&+\displaystyle\frac{a_1}{K}\sqrt{|x|}\displaystyle\sum_{n=0}^{+\infty}\displaystyle\frac{(4K)^n(n!)^2}{(2n)!(2n+1)!}x^{2n+1}.
\end{array}
\end{equation}
Now, for $\lambda=-1$ and $\nu=\displaystyle\frac{3k-k^2-1}{2(k-1)}$ for some $k\in\mathbb{N}$, $k\geq3$ (to guaranty that $\nu\not=\displaystyle\frac{1}{2}$, we get $a_0=a_1=0$ and
\begin{equation}\label{recurence7}
(n-k)(n-\eta)a_n=Ka_{n-2},\;\forall\,n\geq2.
\end{equation}
where $2\eta=1-2\nu-\sqrt{(1-2\nu)^2+4}$. Of course, in the present case $\eta\notin\mathbb{N}$. More precisely, for $\nu<\displaystyle\frac{1}{2}$, we have $0<\eta<1$ and for $\nu>\displaystyle\frac{1}{2}$, it holds that $\eta<0$. Hence, whenever $k=2p$, we obtain $a_{2n+1}=0$ for all $n$ and $a_{2n}=0$ just for $n\leq p-1$. The recurrence relation (\ref{recurence7}) yields that
$$
a_{2n}=a_{2p}\Gamma(n-\zeta)\displaystyle\left(\displaystyle\frac{K}{4}\right)^{n-p}\displaystyle\frac{1}{(n-p)!\Gamma(n-\zeta)},\;n\geq p
$$
where $\zeta=\displaystyle\frac{\eta}{2}$. Consequently,
\begin{equation}\label{hypergeometricsolution7}
f(x)=a_{2p}\Gamma(n-\zeta)|x|^{\nu_k}\displaystyle\sum_{n=p}^{+\infty}\displaystyle\left(\displaystyle\frac{K}{4}\right)^{n-p}
\displaystyle\frac{1}{(n-p)!\Gamma(n-\zeta)}x^{2n},
\end{equation}
where $\nu_k=\displaystyle\frac{(3k-k^2-1)}{2(k-1)}$. Similarly, whenever whenever $k=2p+1$, we obtain $a_{2n}=0$ for all $n$ and $a_{2n+1}=0$ just for $n\leq p-1$ and from the recurrence relation (\ref{recurence7}) we get
$$
a_{2n+1}=a_{2p+1}\displaystyle\frac{\Gamma(2p-\eta+1)}{\Gamma(p-\zeta+1)}K^{n-p}\displaystyle\frac{\Gamma(n-\zeta)}{(n-p)!\Gamma(2n-\eta+1)},\;n\geq p.
$$
Consequently,
\begin{equation}\label{hypergeometricsolution8}
f(x)=a_{2p+1}\displaystyle\frac{\Gamma(2p-\eta+1)}{\Gamma(p-\zeta+1)}|x|^{\nu_k}\displaystyle\sum_{n=p}^{+\infty}
K^{n-p}\displaystyle\frac{\Gamma(n-\zeta)}{(n-p)!\Gamma(2n-\eta+1)}x^{2n+1}.
\end{equation}
It remains now to study the case when $1+2\lambda\nu\not=0$ and $1+2(\lambda+\nu)+2\lambda\nu\not=0$. It holds so that $a_0=a_1=0$ and remains in (\ref{recurence1}) that
\begin{equation}\label{recurence8}
(n^2+(2(\lambda+\nu)-1)n+(1+2\lambda\nu))a_n=Ka_{n-2},\;\forall\,n\geq2.
\end{equation}
So assume now that $(\lambda,\nu)\in D(\Omega,\displaystyle\frac{\sqrt5}{2})$ the open disc of center $\Omega(\displaystyle\frac{1}{2},\displaystyle\frac{1}{2})$ and radius $R=\displaystyle\frac{\sqrt5}{2}$. Then the second order equation $n^2+(2(\lambda+\nu)-1)n+(1+2\lambda\nu))>0$ for all $n$ which yields with the fact that $a_0=a_1=0$ the null solution. \\
Next, assume that $(\lambda,\nu)\in C(\Omega,\displaystyle\frac{\sqrt5}{2})$ the circle of center $\Omega(\displaystyle\frac{1}{2},\displaystyle\frac{1}{2})$ and radius $R=\displaystyle\frac{\sqrt5}{2}$. We get
$$
n^2+(2(\lambda+\nu)-1)n+(1+2\lambda\nu))=(n-\eta)^2,
$$
with $\eta=\displaystyle\frac{1-2(\lambda+\nu)}{2}$. So, whenever there exists $k\in\mathbb{N}$ for which $\eta=2k$, we obtain $a_{2n+1}=0$ for all $n$, and
$$
a_{2n}=a_{2k}\displaystyle\left(\displaystyle\frac{K}{4}\right)^{n-k}\displaystyle\frac{1}{((n-k)!)^2},\;n\geq k
$$
and $a_{2n}=0$ otherwise. Hence, the generated solution will be on the form
\begin{equation}\label{hypergeometricsolution9}
f(x)=a_{2k}|x|^{\nu}\displaystyle\sum_{n=k}^{+\infty}\displaystyle\left(\displaystyle\frac{K}{4}\right)^{n-k}\displaystyle\frac{1}{((n-k)!)^2}x^{2n}
\end{equation}
where $\nu=\displaystyle\frac{1}{2}\pm\sqrt{\displaystyle\frac{5}{4}-(\lambda-\displaystyle\frac{1}{2})^2}$. \\
Now, whenever there exists $k\in\mathbb{N}$ for which $\eta=2k+1$, we obtain $a_{2n}=0$ for all $n$, and
$$
a_{2n+1}=a_{2k+1}(4K)^{n-k}\displaystyle\frac{((n-k)!)^2}{((2n-2k+1)!)^2},\;n\geq k
$$
and $a_{2n}=0$ otherwise. Hence, the generated solution will be on the form
\begin{equation}\label{hypergeometricsolution10}
f(x)=a_{2k+1}|x|^{\nu}\displaystyle\sum_{n=k}^{+\infty}(4K)^{n-k}\displaystyle\frac{((n-k)!)^2}{((2n-2k+1)!)^2}x^{2n+1}
\end{equation}
for the same expression of $\nu=\displaystyle\frac{1}{2}\pm\sqrt{\displaystyle\frac{5}{4}-(\lambda-\displaystyle\frac{1}{2})^2}$ as in the last previous case.\\
Next, assume that $(\lambda,\nu)\notin\overline{D}(\Omega,\displaystyle\frac{\sqrt5}{2})$ the complementary of the closed disc of center $\Omega(\displaystyle\frac{1}{2},\displaystyle\frac{1}{2})$ and radius $R=\displaystyle\frac{\sqrt5}{2}$ and denote
$$
\eta_1=\displaystyle\frac{1-2(\lambda+\nu)+\sqrt{(2\lambda-1)^2+(2\nu-1)^2-5}}{2}
$$
and
$$
\eta_2=\displaystyle\frac{1-2(\lambda+\nu)-\sqrt{(2\lambda-1)^2+(2\nu-1)^2-5}}{2}
$$
so that the recurrence (\ref{recurence1}) becomes
$$
(n-\eta_1)(n-\eta_2)a_n=Ka_{n-2}.
$$
Next, as previously we get one of the following solutions.
$$
a_{2n}=a_{2k}\displaystyle\left(\displaystyle\frac{K}{4}\right)^{n-k}\displaystyle\frac{\Gamma(k-\zeta_2)}{(n-k)!\Gamma(n-\zeta_2)},\;n\geq k
$$
and $a_{2n}=0$ otherwise, where $k=\displaystyle\frac{\eta_1}{2}\in\mathbb{N}$ and $\zeta_2=\displaystyle\frac{\eta_2}{2}$. In this case, of course, $a_{2n+1}=0$ for all $n$. Next, whenever $k=\displaystyle\frac{\eta_1-1}{2}\in\mathbb{N}$, we get
$$
a_{2n+1}=a_{2k+1}\displaystyle\frac{\Gamma(2k-\eta_2+1)}{\Gamma(k-\zeta_2)}\displaystyle\frac{K^{n-k}\Gamma(n-\zeta_2)}{(n-k)!\Gamma(2n-\eta_2+1)},\;n\geq k
$$
Hence, the generated solution are on the respective form
\begin{equation}\label{hypergeometricsolution11}
f(x)=a_{2k}|x|^{\nu}\displaystyle\sum_{n=k}^{+\infty}\displaystyle\left(\displaystyle\frac{K}{4}\right)^{n-k}
\displaystyle\frac{\Gamma(k-\zeta_2)}{(n-k)!\Gamma(n-\zeta_2)}x^{2n}
\end{equation}
and
\begin{equation}\label{hypergeometricsolution12}
f(x)=a_{2k+1}\displaystyle\frac{\Gamma(2k-\eta_2+1)}{\Gamma(k-\zeta_2)}|x|^{\nu}\displaystyle\sum_{n=k}^{+\infty}
\displaystyle\frac{K^{n-k}\Gamma(n-\zeta_2)}{(n-k)!\Gamma(2n-\eta_2+1)}x^{2n+1}
\end{equation}
\subsection{A last hypergeometric in time class of solutions}
In this section we will prove that the EPD system has solutions of the form \begin{equation}\label{last-hypergeometric-in-time-solution}
u(x,y,t)=v(x,y,t)=\psi(t)\varphi(x,y)
\end{equation}
where $\psi$ is an hypergeometric type function. Indeed, substituting in (\ref{LinearContinuousProblem}) we obtain
\begin{equation}\label{ContinuousProblem5}
\left\{\begin{array}{lll}
\psi''(t)+\Gamma_a(t)\psi'(t)-K\psi(t)=0\\
\Delta\varphi+<F_{\lambda,\gamma},\nabla\varphi>-K\varphi=0
\end{array}\right.
\end{equation}
where $K$ is a real constant.\\
The function $\psi$ takes the form of the solution of problem (\ref{ContinuousProblem4}) developed in section 2.2. Next, to develop a solution $\varphi$ of the problem
$$
\Delta\varphi+<F_{\lambda,\gamma},\nabla\varphi>-K\varphi=0
$$
we may proceed as in the previous sections by assuming either
$$
\varphi(x,y)=f(x)+g(y)
$$
or
$$
\varphi(x,y)=f(x)g(y).
$$
In both cases, we get 
$$
\varphi(x,y)=f_H(x)+C,
$$
where $f_H$ is the solution of the homogeneous problem (\ref{ContinuousProblem4}) and $C$ is a suitable constant. Indeed, for the first choice we get the following system  
\begin{equation}\label{ContinuousProblem6}
\left\{\begin{array}{lll}
f''(x)+\Gamma_\lambda(x)f'(x)-Kf(x)=\widetilde{K},\\
g''(y)+\Gamma_\gamma(y)g'(y)-Kg(y)=-\widetilde{K}.
\end{array}\right.
\end{equation}
For $K=0$, we get analogous problem as (\ref{ContinuousProblem3}) by replacing $K$ by $\widetilde{K}$ and thus the solutions may be obtained from section 2.1. Next, whenever $K\not=0$, a particular solution is
$$
(f_{particular},g_{particular})=(-\displaystyle\frac{\widetilde{K}}{K},\displaystyle\frac{\widetilde{K}}{K})
$$
and the homogeneous problem in this case is the same as (\ref{ContinuousProblem4}).
\section{Discrete two-dimensional Nonlinear EPD system}
The object of this section is to explain the discretization scheme proposed to transform problem (\ref{ContinuousProblem1})-(\ref{eqn1-4}) into a discrete
quasi-linear one. Let $n\in\mathbb{N}$ and consider a time step
$l=\Delta t$ and a space one $h=\displaystyle\frac{L_1-L_0}{J+1}$.
Next, denote for $k\in\mathbb{N}$ and $j,m\in\{0,...,J+1\}$
$$
t^k=t_0+kl\;,x_j=L_0+jh\quad\hbox{and}\quad y_m=L_0+mh
$$
so that the cube $[L_0,L_1]\times[L_0,L_1]$ is subdivided into cubes $C_{j,m}=[x_j,x_{j+1}]\times[y_m,y_{m+1}]$. For a function $z$ defined on the cube $[L_0,L_1]\times[L_0,L_1]$, we denote by small $z_{j,m}^k$ the net function $z(x_j,y_m,t^k)$ and capital $Z_{j,m}^k$ the numerical approximation. Consider next the discrete finite difference operators
$$
u_t=\displaystyle\frac{u^{n+1}-u^{n-1}}{2\ell},\;\;u_{tt}=\displaystyle\frac{u^{n+1}-2u^n+u^{n-1}}{\ell^2},
$$
$$
u_x=\displaystyle\frac{\overline{u}^{n}_{j+1,m}-\overline{u}^{n}_{j-1,m}}{2h},
$$
$$
\Delta u
=\displaystyle\frac{\overline{u}^{n}_{j+1,m}-2\overline{u}^{n}_{j,m}+\overline{u}^{n}_{j-1,m}}{h^2}
+\displaystyle\frac{\overline{u}^{n}_{j,m+1}-2\overline{u}^{n}_{j,m}+\overline{u}^{n}_{j,m-1}}{h^2},
$$
$$
\overline{u}^n=\alpha u^{n+1}+(1-2\alpha)u^n+\alpha u^{n-1}.
$$
Denote also $G(u,v)=|u|^{p-1}v$ and $H(u)=|v|^{p-1}u$.

The discretization of the first equation of problem (\ref{ContinuousProblem1}) yields that
$$
\begin{array}{lll}
&&\displaystyle\frac{u^{n+1}_{j,m}-2u^{n}_{j,m}+u_{j,m}^{n-1}}{\ell^2}
+\displaystyle\frac{2a}{t_n}\displaystyle\frac{v^{n+1}_{j,m}-v_{j,m}^{n-1}}{2\ell}\\
&=&\alpha\displaystyle\frac{u^{n+1}_{j+1,m}-2u^{n+1}_{j,m}+u_{j-1,m}^{n+1}}{h^2}+(1-2\alpha) \displaystyle\frac{u^{n}_{j+1,m}-2u^{n}_{j,m}+u_{j-1,m}^{n}}{h^2}\\
&+&\alpha\displaystyle\frac{u^{n-1}_{j+1,m}-2u^{n-1}_{j,m}+u_{j-1,m}^{n-1}}{h^2}
+\alpha\displaystyle\frac{u^{n+1}_{j,m+1}-2u^{n+1}_{j,m}+u_{j,m-1}^{n+1}}{h^2}\\
&+&(1-2\alpha)\displaystyle\frac{u^{n}_{j,m+1}-2u^{n}_{j,m}+u_{j,m-1}^{n}}{h^2}
+\alpha\displaystyle\frac{u^{n-1}_{j,m+1}-2u^{n-1}_{j,m}+u_{j,m-1}^{n-1}}{h^2}\\
&+&\displaystyle\frac{2\lambda}{x_j}[\alpha\displaystyle\frac{v^{n+1}_{j+1,m}-v^{n+1}_{j-1,m}}{2h}+(1-2\alpha) \displaystyle\frac{v^{n}_{j+1,m}-v^{n}_{j-1,m}}{2h}+\alpha\displaystyle\frac{v^{n-1}_{j+1,m}-v^{n-1}_{j-1,m}}{2h}]\\
&+&\displaystyle\frac{2\gamma}{y_m}[\alpha\displaystyle\frac{v^{n+1}_{j,m+1}-v^{n+1}_{j,m-1}}{2h}+(1-2\alpha) \displaystyle\frac{v^{n}_{j,m+1}-v^{n}_{j,m-1}}{2h}+\alpha\displaystyle\frac{v^{n-1}_{j,m+1}-v^{n-1}_{j,m-1}}{2h}]\\
&+&\displaystyle\frac{G_{j,m}^n+G_{j,m}^{n-1}}{2}.
\end{array}
$$
Denote next
$$
{a_n}={\displaystyle\frac{a}{t_n}},\;\;
{\lambda_j}={\displaystyle\frac{\lambda}{x_j}},\;\;
{\gamma_m}={\displaystyle\frac{\gamma}{y_m}},\;\;\sigma=\displaystyle\frac{l^2}{h^2}.
$$
Par suite
$$
\begin{array}{lll}
&&u_{j,m}^{n+1}-2u_{j,m}^{n}+u_{j,m}^{n-1}+\ell{a_n}(v_{j,m}^{n+1}-v_{j,m}^{n-1})\\
&=&\sigma\{\alpha[u^{n+1}_{j+1,m}-2u^{n+1}_{j,m}+u_{j-1,m}^{n+1}]+(1-2\alpha)[u^{n}_{j+1,m}-2u^{n}_{j,m}+u_{j-1,m}^{n}]\\
&+&\alpha[u^{n-1}_{j+1,m}-2u^{n-1}_{j,m}+u_{j-1,m}^{n-1}]+\alpha[u^{n+1}_{j,m+1}-2u^{n+1}_{j,m}+u_{j,m-1}^{n+1}]\\
&+&(1-2\alpha)[u^{n}_{j,m+1}-2u^{n}_{j,m}+u_{j,m-1}^{n}+\alpha[u^{n-1}_{j,m+1}-2u^{n-1}_{j,m}+u_{j,m-14}^{n-1}]]\}\\
&+&\sigma h{\lambda_j}\{\alpha[v^{n+1}_{j+1,m}-v_{j-1,m}^{n+1}]
+(1-2\alpha)[v^{n}_{j+1,m}-v_{j-1,m}^{n}]+\alpha[v^{n-1}_{j+1,m}-v_{j-1,m}^{n-1}]\}\\
&+&\sigma h{\gamma_m}\{\alpha[v^{n+1}_{j,m+1}-v_{j,m-1}^{n+1}]
+(1-2\alpha)[v^{n}_{j,m+1}-v_{j,m-1}^{n}]+\alpha[v^{n-1}_{j,m+1}-v_{j,m-1}^{n-1}]\}
\end{array}
$$
or equivalently,
$$
\begin{array}{lll}
&&u_{j,m}^{n+1}-2 u_{j,m}^{n}+u_{j,m}^{n-1}+\ell a_n (v_{j,m}^{n+1}-v_{j,m}^{n-1})\\
&=&\sigma\{\alpha[u^{n+1}_{j+1,m}-2u^{n+1}_{j,m}+u_{j-1,m}^{n+1}]+(1-2\alpha)[u^{n}_{j+1,m}-2u^{n}_{j,m}+u_{j-1,m}^{n}]\\
&+&\alpha[u^{n-1}_{j+1,m}-2u^{n-1}_{j,m}+u_{j-1,m}^{n-1}]+\alpha[u^{n+1}_{j,m+1}-2u^{n+1}_{j,m}+u_{j,m-1}^{n+1}]\\
&+&(1-2\alpha)[u^{n}_{j,m+1}-2u^{n}_{j,m}+u_{j,m-1}^{n}+\alpha[u^{n-1}_{j,m+1}-2u^{n-1}_{j,m}+u_{j,m-1}^{n-1}]]\}\\
&+&\sigma h\lambda_j\{\alpha[v^{n+1}_{j+1,m}-v_{j-1,m}^{n+1}]+(1-2\alpha)[v^{n}_{j+1,m}-v_{j-1,m}^{n}]+\alpha[v^{n-1}_{j+1,m}-v_{j-1,m}^{n-1}]\}\\
&+&\sigma h\gamma_m\{\alpha[v^{n+1}_{j,m+1}-v_{j,m-1}^{n+1}]+(1-2\alpha)[v^{n}_{j,m+1}-v_{j,m-1}^{n}]+\alpha[v^{n-1}_{j,m+1}-v_{j,m-1}^{n-1}]\}
\end{array}
$$
This may be written in a matrix-vector form
$$
\begin{array}{lll}
&&IU^{n+1}-2IU^{n}+IU^{n-1}+a_n I(V^{n+1}-V^{n-1})\\
&=&\sigma\{\alpha AU^{n+1}+(1-2\alpha)AU^n+\alpha AU^{n-1}+\alpha U^{n+1}A+(1-2\alpha)U^nA+\alpha U^{n-1}A\}\\
&+&\Theta\{\alpha V^{n+1}+(1-2\alpha)V^{n}+\alpha V^{n-1}\}+\{\alpha V^{n+1}+(1-2\alpha)V^{n}+\alpha V^{n-1}\}\Lambda
\end{array}
$$
where $A$, $\Theta$ and $\Lambda$ are the matrices given by
$$
A_{1,2}=A_{J+1,J}=2,\quad\Theta_{1,2}=\Theta_{J+1,J}=\Lambda_{1,2}=\Lambda_{J+1,J}=0,
$$
and for $1\leq j,m\leq J$, 
$$
A_{jj}=-2,\quad A_{j,j+1}=A_{j,j-1}=1\quad\hbox{and}\quad0\;\hbox{elsewhere},
$$
$$
\Theta_{j,j}=0,\quad \Theta_{j,j+1}=-\Theta_{j,j-1}=\lambda_j\quad\hbox{and}\quad0\;\hbox{elsewhere},
$$
and
$$
\Lambda_{j,j}=0,\quad \Lambda_{j+1,j}=-\Lambda_{j-1,j}=\delta_j\quad\hbox{and}\;0\;\hbox{elsewhere},
$$
Denote next,
$$
{W_{\alpha}}={\left(\displaystyle\frac{1}{2}I-\alpha\sigma A\right)}
$$
$$
{R_{n,\alpha}}={\left(\displaystyle\frac{1}{2}a_nI-\alpha\sigma h\Theta\right)}
$$
$$
{S_{n,\alpha}}={\left(\displaystyle\frac{1}{2}a_nI-\alpha\sigma h\Lambda\right)}.
$$
We get
\begin{equation}\label{LyapEq1}
\begin{array}{lll}
&&W_{\alpha}U^{n+1}+U^{n+1}W_{\alpha}+R_{n,\alpha}V^{n+1}+V^{n+1}S_{n,\alpha}\\
&=&2\left[W_{\alpha-\frac{1}{2}}U^n+U^nW_{\alpha-\frac{1}{2}}\right]-\left[W_\alpha U^{n-1}+U^{n-1}W_\alpha\right]\\
&&+(1-2\alpha)\left[\Theta V^n+V^n\Lambda\right]+R_{n,-\alpha}V^{n-1}+V^{n-1}S_{n,-\alpha}\\
&&+\displaystyle\frac{G^n+G^{n-1}}{2}.
\end{array}
\end{equation}
Similarly, for the second equation we get
\begin{equation}\label{LyapEq2}
\begin{array}{lll}
&&W_{\alpha}V^{n+1}+V^{n+1}W_{\alpha}+R_{n,\alpha}U^{n+1}+U^{n+1}S_{n,\alpha}\\
&=&2\left[W_{\alpha-\frac{1}{2}}V^n+V^nW_{\alpha-\frac{1}{2}}\right]-\left[W_\alpha V^{n-1}+V^{n-1}W_\alpha\right]\\
&&+(1-2\alpha)\left[\Theta U^n+U^n\Lambda\right]+R_{n,-\alpha}U^{n-1}+U^{n-1}S_{n,-\alpha}\\
&&+\displaystyle\frac{H^n+H^{n-1}}{2}.
\end{array}
\end{equation}
\section{Solvability of the discrete problem}
In \cite{Benmabrouk1}, the authors have transformed the Lyapunov operator
obtained from the discretization method into a standard linear operator
acting on one column vector by juxtaposing the columns of the matrix $X$
horizontally which leads to an equivalent linear operator characterized by a
fringe-tridiagonal matrix. We used standard computation to prove the
invertibility of such an operator. Here. we do not apply the same
computations as in \cite{Benmabrouk1}, but we develop different arguments.
The first main result is stated as follows.

\begin{theorem}\label{theorem1}
The system (\ref{LyapEq1})-(\ref{LyapEq2}) is uniquely solvable whenever $U^0$ and $U^1$ are known.
\end{theorem}
\textbf{Proof.} It reposes on the inverse of Lyapunov-Syslvester operators. Consider the endomorphism $\Phi$ defined by 
\begin{equation}\label{LyapSylOperator}
\Phi_{l,h}(X,Y)=(\mathcal{L}_{W_\alpha}(X)+\mathcal{L}_{R_{n,\alpha},S_{n,\alpha}}(Y,\mathcal{L}_{W_\alpha}(Y)+\mathcal{L}_{R_{n,\alpha},S_{n,\alpha}}(X)),
\end{equation}
where for matrices $A$ and $B$, $\mathcal{L}_{A,B}$ is the Lyapunov-Sylvester operator defined by 
$$
\mathcal{L}_{A,B}(X)=AX+XB
$$
and where $\mathcal{L}_{A}\equiv\mathcal{L}_{A,A}$. To prove Theorem \ref{theorem1}, it suffices to show that $ker\Phi $ is reduced to $0$. Indeed, whenever $l=o(h)$ and $l,h\rightarrow0$, we get
$$
\Phi_{l,h}(X,Y)\rightarrow\Phi(X,Y)=(X+a_nY,Y+a_nX).
$$
Next, whenever $\Phi(X,Y)=0$, we get
$$
X+a_nY=Y+a_nX=0.
$$
Which means that $X=Y=0$. Next, we apply the following result.
\begin{lemma} \label{LemmeInversion}
Let $E$ be a finite dimensional ($\mathbb{R}$ or $\mathbb{C}$) vector space and $(\Phi_n)_n$ be a sequence of endomorphisms converging uniformly to an invertible endomorphism $\Phi$. Then, there exists $n_{0}$ such that, for any $n\geq\,n_{0}$, the endomorphism $\Phi_n$ is invertible.
\end{lemma}
Observing that the operator $\Phi$ obtained above is invertible, we get that $\Phi_{l,h}$ is invertible for $l,h$ small enough.
\section{Consistency, stability and convergence of the discrete method}
Recall firstly that the consistency of the numerical scheme is always done by evaluating the local truncation error arising from the discrete and the continuous problem. Applying Taylor's expansion in the discrete equations raised in section 2, we get the following truncation principal part for the first equation in system (\ref{ContinuousProblem1})
$$
\begin{array}{lll}
\mathcal{L}_{u,v}^1(x,y,t)
&=&l^2\Bigl[\displaystyle\frac{1}{12}\frac{\partial^4u}{\partial t^4}+\displaystyle\frac{1}{6}\Gamma_a(t)\displaystyle\frac{\partial^3u}{\partial t^3}+\alpha\frac{\partial^2}{\partial t^2}(\Delta u+2<F_{\lambda,\gamma,\nabla v}>)\Bigr]\\
\\
&&+h^2\Bigl[\displaystyle\frac{1}{12}\Delta_2u+\displaystyle\frac{1}{6}<F_{\lambda,\gamma,\nabla_3v}>\Bigr]+o(l^2+h^2)
\end{array}
$$
and for the second equation, we get
$$
\begin{array}{lll}
\mathcal{L}_{u,v}^2(x,y,t)
&=&l^2\Bigl[\displaystyle\frac{1}{12}\frac{\partial^4v}{\partial t^4}+\displaystyle\frac{1}{6}\Gamma_a(t)\displaystyle\frac{\partial^3v}{\partial t^3}+\alpha\frac{\partial^2}{\partial t^2}(\Delta v+2<F_{\lambda,\gamma,\nabla u}>)\Bigr]\\
\\
&&+h^2\Bigl[\displaystyle\frac{1}{12}\Delta_2v+\displaystyle\frac{1}{6}<F_{\lambda,\gamma,\nabla_3u}>\Bigr]+o(l^2+h^2),
\end{array}
$$
where $\Delta_2=\displaystyle\frac{\partial^4}{\partial x^4}+\frac{\partial^4}{\partial y^4}$ and $\nabla_3=(\displaystyle\frac{\partial^3}{\partial x^3},\frac{\partial^3}{\partial y^3})$.
WE then have the following lemma.
\begin{lemma}
	The discrete scheme is consistent with order $(l^2+h^2)$.
\end{lemma}
Now, we will examine the stability of the scheme. We will apply the same method as in \cite{Benmabrouk1}, \cite{Bezia-BenMabrouk-Betina1}, \cite{Bezia-BenMabrouk1} and \cite{Chteouietal1} based on the Lyapunov criterion of stability. A linear system $\mathcal{L}(u_{n+1},u_{n},u_{n-1},\dots)=0$ is stable in the sense of Lyapunov iff for any bounded initial value $u_0$, the solution $u_n$ ramains bounded for all $n\geq0$.
\begin{lemma} \label{LyapunovStabilityLemma}
	$\mathcal{P}_n$: The solution $(U^n,V^n) $ is bounded independently of $n$ whenever the initial solution $(U^0,V^0)$ is bounded.
\end{lemma}
Before going on proving this result, we stress on the fact that contrarily to previous studies such as \cite{Benmabrouk1}, \cite{Bezia-BenMabrouk-Betina1}, \cite{Bezia-BenMabrouk1} and \cite{Chteouietal1}, we are confronted in the present study to Lyapunov-Sylvester operators that are non commutative, which raises some new difficulties in the proof of the stability and thus imposes different ideas.\\
\textbf{Proof.} Recall firstly that the discrete scheme in the matrix form may be written as
$$
\Phi_{l,h}(U^{n+1},V^{n+1})=F(U^n,V^n,U^{n-1},V^{n-1}),
$$
where $\Phi_{l,h}$ is the Lyapunov-Sylvester operator defined in (\ref{LyapSylOperator}) and where the linear operator $F=(F_1,F_2)$ (which is a melange of Lyapunov-Sylvester operators) is defined by
$$
\begin{array}{lll}
\medskip F_1(X,Y,Z,T)&=&2\mathcal{L}_{W_{\alpha-1/2}}(X)-\mathcal{L}_{W_{\alpha}}(Z)\hfill\crcr\medskip
&&+(1-2\alpha)\mathcal{L}_{\Theta,\Gamma}(Y)
+\mathcal{L}_{R_{n,-\alpha},S_{n,-\alpha}}(T)\hfill\crcr\medskip
&&+G(X,Y,Z,T)
\end{array}
$$
and
$$
\begin{array}{lll}
\hfill\crcr\medskip
F_2(X,Y,Z,T)&=&2\mathcal{L}_{W_{\alpha-1/2}}(Y)-\mathcal{L}_{W_{\alpha}}(T)\hfill\crcr\medskip
&&+(1-2\alpha)\mathcal{L}_{\Theta,\Gamma}(X)
+\mathcal{L}_{R_{n,-\alpha},S_{n,-\alpha}}(Z)\hfill\crcr\medskip
&&+H(X,Y,Z,T),
\end{array}
$$
where
$G(X,Y,Z,T)$ is the matrix with coefficients 
$$
G(X,Y,Z,T)_{j,m}=\displaystyle\frac{|X_{j,m}|^{p-1}Y_{j,m}+|Z_{j,m}|^{p-1}T_{j,m}}{2}
$$
and $H(X,Y,Z,T)$ is the matrix with coefficients 
$$
H(X,Y,Z,T)_{j,m}=\displaystyle\frac{|Y_{j,m}|^{q-1}X_{j,m}+|T_{j,m}|^{q-1}Z_{j,m}}{2}.
$$
We remark immediately that
$$
\|G(X,Y,Z,T)\|\leq\displaystyle\frac{\|X\|^{p-1}\|Y\|+\|Z\|^{p-1}\|T\|}{2}
$$
and  
$$
\|H(X,Y,Z,T)\|\leq\displaystyle\frac{\|Y\|^{q-1}\|X\|+\|T\|^{q-1}\|Z\|}{2}.
$$
Consequently,
$$
\|G(X,Y,Z,T)\|\leq\displaystyle\frac{\|(X,Y)\|^{p}+\|(Z,T)\|^{p}}{2}
\leq\displaystyle\|(X,Y,Z,T)\|^{p}
$$
and  
$$
\|H(X,Y,Z,T)\|\leq\displaystyle\frac{\|(X,Y)\|^{q}+\|(Z,T)\|^{q}}{2}
\leq\displaystyle\|(X,Y,Z,T)\|^{q}.
$$
Recall also that
$$
\Phi_{l,h}(X,Y)\rightarrow\Phi(X,Y)=(X+a_nY,Y+a_nX),\quad\hbox{as}\quad l,h\rightarrow0.
$$
Observe next that for $l,h$ small enough
$$
\|\Phi\|\geq|a_n-1|,\;\forall\,n.
$$
So, for $n$ large enough, we may obtain
$$
\|\Phi\|\geq\displaystyle\frac{1}{2}.
$$
On the other hand, we have
$$
\|\Phi_{l,h}(X,Y)-\Phi(X,Y)\|\geq\alpha\sigma\left(2\|A\|+h\|\Theta\|+h\|\Gamma\|\right)\|X+Y\|.
$$
Next, observe that 
$$
(2\|A\|+h\|\Theta\|+h\|\Gamma\|)\leq4+h(\displaystyle\max_j|\lambda_j|+\displaystyle\max_j|\gamma_j|)
$$
and consequently, by denoting 
$$
C_\alpha=\alpha[4+h(\displaystyle\max_j|\lambda_j|+\displaystyle\max_j|\gamma_j|)],
$$
we get
$$
\|\Phi_{l,h}(X,Y)-\Phi(X,Y)\|\geq C_\alpha\|(X,Y)\|.
$$
Conseqyently,
$$
\left(\displaystyle\frac{1}{2}-C_\alpha\sigma\right)\|(X,Y)\|\leq\|\Phi_{l,h}(X,Y)\|.
$$
Next, by choosing the time step $l$ and the space step $h$ so that $4\sigma C_\alpha<1$, we obtain
$$
\displaystyle\frac{1}{4}\|(X,Y)\|\leq\|\Phi_{l,h}(X,Y)\|.
$$
Finally, the lemma follows by recurrence on $n$ by observing for the choice above that
$$
\|(U^{n+1},V^{n+1})\|\leq4\|F(U^n,V^n,U^{n-1},V^{n-1})\|
$$
and the continuity of $F$.

Now, it remains finally to check the convergence of the discrete scheme. This is done by a direct application of the following well-known result \cite{Lax-Richtmyer}.
\begin{theorem} 
(\textbf{Lax Equivalence Theorem}). For a consistent finite difference scheme, stability is equivalent to convergence.
\end{theorem}
\begin{lemma} \label{laxequivresult}
As the numerical scheme is consistent and stable, it is then convergent.
\end{lemma}
\section{Numerical implementation}
In this section we present a numerical example that links between the EPD system studied here and the classical EPD equation in some parts so that we point out two aims; we firstly validate the theoretical results developed previously and secondly to show that such a system may be considered as a perturbation of the classical EPD equation. 

To measure the closeness of the numerical solution and the exact one, the error is evaluated via an $L_2$ matrix norm 
$$
\|X\|_2=\Big(\sum_{i,j=1}^{J+2}|X_{ij}|^2\Big)^{1/2}
$$
for a matrix $X=(X_{ij})\in\mathcal{M}_{J+2}\mathbb{C}$. Denote $u^n$ the net function $u(x,y,t^n)$ and $U^n$ the numerical solution. We propose to compute the discrete error
\begin{equation}
\mathrm{Er}=\max_n\|U^n-u^n\|_2  \label{Er}
\end{equation}
on the grid $(x_{i},y_{j})$, $0\leq\,i,j\leq J+1$ and the relative error between the exact solution and the numerical one as
\begin{equation}
\mathrm{Relative\,Er}=\max_n\frac{\|U^n-u^n\|_2}{\|u^n\|_2}  \label{Errelative}
\end{equation}
on the same grid. Next, we will assume that 
$$
\lambda=\gamma\displaystyle\frac{1}{4},\;a=\displaystyle\frac{5}{2},\;
p=\displaystyle\frac{3}{2}\;\;\hbox{and}\;\;q=\displaystyle\frac{4}{3}.
$$
We denote also
$$
r^2=x^2+y^2\quad\hbox{and}\quad\,g_p(x,y,t)=\exp\left(-p(\displaystyle\frac{t^2}{2}+r^2)\right).
$$
We consider the inhomogeneous system 
\begin{equation}\label{InhpmogeneousContinuousProblem1}
\left\{\begin{array}{lll}
u_{tt}+\Gamma_a(t)v_t=\Delta u+<F_{\lambda,\gamma},\nabla v>+|u|^{p-1}v+G_1(x,y,t)\\
v_{tt}+\Gamma_a(t)u_t=\Delta v+<F_{\lambda,\gamma},\nabla u>+|v|^{q-1}u+G_2(x,y,t).
\end{array}\right.
\end{equation}
where $G_1$ and $G_2$ are explicited respectively by
$$
G_1(x,y,t)=\left(t^2-4r^2\right)g_1(x,y,t)-g_p(x,y,t)
$$
and
$$
G_2(x,y,t)=\left(t^2-4r^2\right)g_1(x,y,t)-g_q(x,y,t).
$$
The exact solution of such a system is given by 
$$
u(x,y,t)=v(x,y,t)=
\exp\left(-(\displaystyle\frac{t^2}{2}+r^2)\right).
$$
We fix the domain bounds as follows:
$$
t_0=0,\quad\mbox{and}\quad L_1=-L_0=10.
$$
Finally tu fulful the assumption $l=o(h)$ we assume that $l=h\sqrt{h}$. The following table resumes the error estimates, the relative error as well as the time of execusion for the corresponding algorithm (denoted here II) compared to the time execusion of the classical tri-diagonal one (denoted I) for different values of the space step.
\begin{table}[th]
\caption{} \label{table1}
\begin{center}
\begin{tabular}{|c|c|c|c|c|c|c|}
\hline
J & ErI & ErII & Relative ErI & Relative ErII & TimeI & TimeII\\ 
\hline
24 & 4.01E-3 & 3.31E-3 & 0.1415 & 0.1382 & 20\,s & 05\,s \\ \hline
49 & 4.01E-3 & 3.10E-4 & 0.1423 & 0.1380 & 49\,s & 7.1\,s \\ \hline
99 & 3.80E-3 & 3.20E-5 & 0.1431 & 0.1381 & 466\,s & 20\,s \\ \hline
149 & 3.10E-4 & 3.03E-5 & 0.1433 & 0.1383 & 1002\,s & 31.2\,s \\ \hline
199 & 2.80E-4 & 2.19E-5 & 0.1436 & 0.1384 & 2114\,s & 35.6\,s \\ \hline
399 & 2.30E-4 & 1.29E-6 & 0.1441 & 0.1384 & 3215\,s & 41.2\,s \\ \hline
499 & 2.10E-4 & 1.11E-6 & 0.1445 & 0.1385 & 4001\,s & 46\,s \\ \hline
999 & 2.01E-4 & 1.02E-6 & 0.1446 & 0.1386 & 4516\,s & 49.3\,s \\ \hline
\end{tabular}
\end{center}
\end{table}
In numerical studies of PDEs one important task is the convergence of algorithms and especially the rate of convergence. Different methods have been developed to get fast algorithms. In the present work, Table 1 above includes a comparison between existing method (denote Method I) and the present scheme based on the generalized Lyapunov-Syslvester form (denote method II) for different values of the parameter $J$. It is remarkable that our method is faster than the classical one. This is perfect as nowadays focuses are on big and/or cloud data and thus may seek fast and accurate algirithms. In Table 1, we noticed easily that an accelerated procedure is pointed out. For $J=24$ for example, a 4-times faster algorithm is obtained. For $J=99$ the running time is reduced to be more than 20-times. Increasing more the mesh size ($J$) results in more and more best faster algorithms with a rate of running time proportion crossing 90-times for $J=999$. 
\section{Conclusion}
In this work, computational method has been developed for numerical solutions of Poisson-Darbous-Euler system in 2-D case. The method has yielded algebraic systems based on general Lyapunov-Sylvester equations. Fast and accurate algorithms have been obtained compared with the associated tri-diagonal classical algorithms always applied in such problems. 

\end{document}